\documentclass{amsart}
\usepackage{amsfonts,yfonts,amsmath,amsthm,amssymb,bbm}
\usepackage[all]{xy}
\newtheorem{theorem}{Theorem}
\newtheorem{lemma}{Lemma}
\newtheorem{proposition}{Proposition}

\newcommand{\rl}{{\mathbb{R}}}
\newcommand{\cx}{{\mathbb{C}}}
\newcommand{\id}{{\mathbb{I}}}
\newcommand{\m}{{\mathcal{M}}}
\newcommand{\dbar}{\overline{\partial}}
\newcommand{\Db}[1]{\frac{\partial{#1}}{\partial\overline{z}}}

\newcommand{\assign}{:=}
\newcommand{\abs}[1]{\left|{#1}\right|}
\newcommand{\norm}[1]{\left\|{#1}\right\|}
\newcommand{\cts}[2]{\mathcal{A}\left({#1},{#2}\right)}
\newcommand{\ctss}[2]{\mathcal{A}_\phi\left({#1},{#2}\right)}
\newcommand{\hol}[2]{\mathcal{O}\left({#1},{#2}\right)}
\newcommand{\hols}[2]{\mathcal{O}_\phi\left({#1},{#2}\right)}
\newcommand{\tmop}[1]{\ensuremath{\operatorname{#1}}}
\newcommand{\dist} {\tmop{dist}}

\newcommand{\intr}[1]{{#1}^{\circ}}
\newenvironment{Proof}{\begin{proof}
  }{\end{proof}}
\newenvironment{myprop}[1]{\mbox{}\newline{\bf{#1}.}\em}{\mbox{}\newline}

\newcounter{romct}

\newcommand{\aref}[1]{\S\ref{#1}}
\title[Sets of Approximation and Interpolation ]
{Sets of Approximation and Interpolation in $\cx$  for
manifold-valued maps}
\author{Debraj Chakrabarti}
\begin{document}
\maketitle
\bibliographystyle{plain}

\section{Introduction}
\label{intro} \subsection{Notation} This article is devoted to the
study of approximation and interpolation of maps from a compact
$K\subset\cx$ into a complex manifold $\m$, and in particular to
giving examples of sets $K$ for which certain types of approximation
and interpolation are possible. For brevity, we introduce three
properties $A_1$, $A_2$ and $A_3$ that compact subsets of the plane
may possess. We first define property $A_2$.

For a compact $K\subset \cx$, let $\intr{K}$ be the interior of $K$.
We will let $\cts{K}{\m}$ denote the continuous maps from $K$ into
the complex manifold $\m$ which are holomorphic on $\intr{K}$. Let
$\hol{K}{\m}\subset\cts{K}{\m}$ denote the subspace of those maps
$f$ which extend to a holomorphic map from some neighborhood $U_f$
of $K$ to $\m$. We endow $\m$ with an arbitrary metric, which makes
$\cts{K}{\m}$ and $\hol{K}{\m}$ into metric spaces with the uniform
metric on maps.  We will say that a compact $K\subset\cx$ has the
property $A_2$, if, {\em for every complex manifold $\m$, every
finite set $\mathcal{P}\subset\m$, and every $\epsilon>0$, we can
approximate any $f\in\cts{K}{\m}$ by a map
$f_\epsilon\in\hol{K}{\m}$ such that
$\dist_\m(f,f_\epsilon)<\epsilon$, and for each $p\in\mathcal{P}$,
we have $f_\epsilon(p)=f(p)$.} That is, $f_\epsilon$ is an uniform
approximation to $f$ which interpolates the values of $f$ on
$\mathcal{P}$.

We now define the property $A_1$, which will be the special case of
$A_2$ with $\m=\cx$. More precisely, we will say that $K\subset\cx$
has the property $A_2$, if, {\em for any finite $\mathcal{P}\subset
K$, any function in $\mathcal{A}(K)\assign \cts{K}{\cx}$ can be
uniformly approximated by  functions in
$\mathcal{O}(K)\assign\hol{K}{\cx}$ which interpolate the values of
$f$ on $\mathcal{P}$.} Obviously for a set $K$,  $A_2\Rightarrow
A_1$.

It is easy to see that a set $K$ has property $A_1$ iff
$\mathcal{O}(K)$ is dense in $\mathcal{A}(K)$. Note the trivial fact
that there is a constant $C>0$ such that for $g\in\mathcal{C}(K)$,
if $L_\mathcal{P}(g)$ is the Lagrange polynomial which interpolates
the values of $g$ on $\mathcal{P}$, we have
$\norm{L_\mathcal{P}(g)}_K<C\norm{g}_K$. If
$\tilde{f}\in\mathcal{O}(K)$ be such that
$\abs{f-\tilde{f}}<\frac{\epsilon}{C+1}$, then
$\abs{f-f_\epsilon}<\epsilon$, and $f(p)=f_\epsilon(p)$ for
$p\in\mathcal{P}$, where $f_\epsilon= \tilde{f} +
L_\mathcal{P}(f-\tilde{f})$.

The compact sets $K\subset\cx$ such that $\mathcal{O}(K)$ is dense
in $\mathcal{A}(K)$ can be characterized by Vitu\v{s}kin's theorem (
\cite{vitushkin}). A sufficient condition  is that $\cx\setminus K$
has finitely many connected components. Also, such approximation can
be localized, in the sense that $\mathcal{O}(K)$ is dense in
$\mathcal{A}(K)$  iff every point $z\in K$ has a neighborhood $U_z$
in $\cx$ such that $\mathcal{O}(K\cap \overline{U_z})$ is dense in
$\mathcal{A}(K\cap \overline{U_z})$.

We now define property $A_3$. Let $\phi:\m\rightarrow\cx$ be a
holomorphic submersion such that $\phi(K)\supset \m$. Let
$\ctss{K}{\m}$ (resp. $\hols{K}{\m}$ ) be the subspace of
$\cts{K}{\m}$ (resp. $\hol{K}{\m}$) consisting of sections of $\phi$
over $K$, i.e. maps $s:K\rightarrow\m$ such that $\phi\circ s=
\id_K$. We will say that a compact $K\subset \cx$ has property $A_3$
if {\em for every complex manifold $\m$ and every holomorphic
submersion $\phi$ such that $\phi(\m)\supset K$, and every finite
set $\mathcal{P}\subset K$, every section $\sigma\in\ctss{K}{\m}$
can be uniformly approximated by sections in $\hols{K}{\m}$ which
interpolate the values of $\sigma$ on $\mathcal{P}$}.

We also note the following elementary fact:
\begin{lemma}\label{aone} For a compact $K\subset \cx$,
$A_3\Rightarrow A_2 $.\end{lemma}
\begin{Proof}
Assume that  $K$ satisfies $A_3$, and let $\mathcal{P}$ and $\m$
have the same meaning as above, and let $f\in\cts{K}{\m}$. Consider
the complex manifold
 $\mathcal{N}=\m\times \cx$ and let $\phi$ and $\pi$ be the
 projections onto $\cx$ and $\m$ respectively. If $F:K\rightarrow
 \mathcal{N}$ is defined by $F(z)=(f(z),z)$, then clearly
 $F\in\ctss{K}{\mathcal{N}}$ and since $K$ has the property $A_3$,
 we can approximate $F$ by maps $G\in \hols{K}{\mathcal{N}}$ such
 that $G(p)=F(p)=(f(p),p)$ for each $p\in\mathcal{P}$. Then
 $g=\pi\circ G$ is in $\hol{K}{\m}$ and an approximation to $f$ with
 $g(p)=f(p)$ for each $p\in\mathcal{P}$,
  which shows that $K$ has property $A_2$.
\end{Proof}
(We may remark here that while holomorphic submersions
$\phi:\m\rightarrow\cx$ always exist if $\m$ is Stein
(\cite{forstneric:noncritical}), in general such $\m$ may be highly
nontrivial, see e.g. \cite{demailly:nonsteinbundle}.)

It is natural to ask the question: {\em which are the compact  sets
in the plane which satisfy property $A_2$ or property $A_3$?} Given
that approximation is local in nature, one can conjecture that using
patching arguments in the manifold $\m$ we may be able to  show that
every set with property $A_1$ has property $A_2$, or even $A_3$. We
do not know if this program can be carried out. In this article, we
confine ourselves to to the much more modest goal of giving examples
of sets $K$ for which properties $A_2$ and $A_3$ hold (see
Theorems~\ref{onedim},\ref{atwo} and \ref{athree} below.)

\subsection{Known results}
Some  results have been obtained recently regarding the
approximation and interpolation of manifold-valued maps.

The following  was proved in \cite{michiganpaper}. Define a {\em
Jordan domain} to be an open set $\Omega\Subset\cx$, such that the
boundary $\partial\Omega$ has finitely many connected components,
each of which is homeomorphic to a circle. (We explicitly allow
$\Omega$ to be multiply connected.)
\begin{proposition}\label{thesis} Let $K=\overline{\Omega}$, where
 $\Omega\Subset\cx$ is a Jordan domain. For any
complex manifold $\m$, the subspace $\hol{K}{\m}$ is dense in
$\cts{K}{\m}$.
\end{proposition}
This is purely a statement about approximation, no interpolation
being involved. The following result is contained in what was proved
by Drinovec-Drnov\v{s}ek and Forstneri\v{c} in
\cite{bddforstneric:holocurves}, Theorem 5.1:
\begin{proposition}\label{forst}
Let $\Omega\Subset\cx$ have $\mathcal{C}^2$ boundary, $\m$ be a
complex manifold, and let $f:\overline{\Omega}\rightarrow\m$ be a
map of class $\mathcal{C}^r$ ($r\geq 2$) which is holomorphic in
$\Omega$. Given finitely many points $z_1,\ldots, z_l\in\Omega$, and
an integer $k\in\mathbb{}{N}$, there is a sequence of maps
$f_\nu\in\hol{K}{\m}$ such that $f_\nu$ agrees with $f$ to order $k$
at $z_j$ for $j=1,\ldots,l$ and $\nu\in\mathbb{N}$, and the sequence
$f_\nu$ converges to $f$ in $\mathcal{C}^r(\overline{\Omega})$ as
$\nu\rightarrow\infty$.
\end{proposition}
This can in fact be proved when $\Omega\Subset S$, where $S$ is a
non-compact Riemann surface, and $\m$ is a complex space. Observe
that there is a stronger assumption on the smoothness of the map
than in property $A_2$, and also a stronger conclusion regarding
approximation and interpolation. However, the set of points of
interpolation is restricted to the interior of $\overline{\Omega}$.
The authors subsequently proved the following:
\begin{proposition}\label{forst2}
If $\Omega\Subset\cx$ is a domain  with  $\mathcal{C}^2$ boundary,
$\m$ a complex manifold, and $\phi:\m\rightarrow\cx$ a holomorphic
submersion such that $\phi(\m)\supset\overline{\Omega}$, then
$\hols{\overline{\Omega}}{\m}$ is dense in
$\ctss{\overline{\Omega}}{\m}$.
\end{proposition}
This is a (very) special case of \cite{bddforstneric:approximation},
Theorem 5.1 regarding the approximation of sections of submersions
over strongly pseudoconvex sets in Stein manifolds, which may be
thought of as a far-reaching generalization of the
Henkin-Ram\'{\i}rez-Kerzman approximation theorem for functions on
such domains (\cite{henkinleiterer:book}, Theorem 2.9.2. )
\subsection{Sets with property $A_2$ or $A_3$} We recall some definitions
from elementary point-set topology. Let $X$ be a topological space.
For an integer $n\geq 0$, we say that $\tmop{dim}(X)\leq n$ if the
following holds : given any open cover $\mathcal{A}$ of $X$, there
is a refinement $\mathcal{B}$ of $\mathcal{A}$ such that any point
of $X$ is contained in at most $n+1$ of the sets in the cover
$\mathcal{B}$. If $\tmop{dim}(X)\leq n$ holds, but
$\tmop{dim}(X)\leq n-1$ does not, we say that the  dimension
$\tmop{dim}(X)$ of $X$ is $n$.  The dimension is clearly a
topological invariant of $X$. Also, if $\tmop{dim}(X)=n$, then each
connected component of $X$  (or more generally, any closed subset)
has dimension less than or equal to $n$.

We say that a topological space is {\em locally contractible} if
there is a basis of the topology consisting of contractible open
sets. We can now state the following:
\begin{theorem}\label{onedim}
Let $K\subset \cx$ be a connected compact set, such that
$\tmop{dim}(K)=1$, and suppose that $K$ is locally contractible.
Then $K$ has property $A_3$.
\end{theorem}
Examples of sets that satisfy the hypotheses are
\begin{enumerate}
\item arcs in $\cx$,  where By an {\em arc} $\alpha$ in a topological
space $X$ we mean a continuous {\em injective} map
$\alpha:[0,1]\rightarrow X$ from the unit interval. The injectivity
avoids pathologies like space filling curves. By standard abuse of
language, we will refer to the image $\alpha([0,1])$ as the arc
$\alpha$.
\item planar realizations of connected {\em graphs} i.e. finite simplicial complexes
having only $0$-simplices (``vertices") and $1$-simplices
(``edges".) We can think of these as unions of arcs in the plane,
any two of which meet at at most one endpoint of each.
\end{enumerate}
Let $K\subset \cx$, then $\tmop{dim}(K)\leq 1$ iff
$\intr{K}=\emptyset$ (\cite{ems:gtopology1}, p.133,Theorem 20.)
However, it is possible for a set $K$ with $\intr{K}=\emptyset$,  to
be {\em not} locally contractible, and yet to satisfy $A_1$ (recall
that $A_1$ is equivalent to $\mathcal{O}(K)$ being dense in
$\mathcal{A}(K)$.) An example is the ``bouquet of circles"
$K=\bigcup_{n=1}^{\infty} C_n$, where $C_n$ is the circle in the
plane with center at $\frac{1}{n}+i0$, and radius $\frac{1}{n}$. It
can be shown that this $K$ has property $A_1$, yet
Theorem~\ref{onedim} does not apply to it. (We don't know if $K$ has
property $A_3$ or $A_2$.)

For completeness, we record  the following easy fact:
\begin{lemma}\label{zerodim}
Let $K$ be a  compact subset of $\cx$ with $\tmop{dim}(K)=0$. Then
$K$ has property $A_3$.
\end{lemma}

We now go on to give examples of two-dimensional sets which have
properties $A_2$ or $A_3$.  Let $k\geq 0$ be an integer, $\infty$ or
$\omega$. We introduce a class of sets in the plane denoted by
$\mathfrak{C}_k$. We will let $\mathfrak{C}_k$ denote the class of
compact sets $K\subset\cx$ such that there is an integer $N$ and
domains $\Omega_j$ , where $j=1,\ldots, N$ with the following
properties:
\begin{itemize}
\item For $k\not=0$, each $\Omega_j$ is a domain with $\mathcal{C}^k$
boundary, and for $k=0$, each $\Omega_j$ is a  Jordan domain.
\item  The $\Omega_j$'s are pairwise disjoint :
$\Omega_i\cap \Omega_j=\emptyset$ if $i\not=j$.
\item $K=\cup_{i=1}^N \overline{\Omega_j}$. We will refer to each
$\Omega_j$ as a {\em summand} of $K$.
 \item whenever $i\not = j$, the set
 \[P_{ij}\assign\overline{\Omega_i}\cap \overline{\Omega_j}=\partial\Omega_i\cap \partial\Omega_j \]
  is finite.
 \item If $k\not= 0$, at each point of $P_{ij}$, the boundaries
  $\partial\Omega_i$ and $\partial\Omega_j$ are tangent to each other.
\end{itemize}
 We can now provide a supply of sets which have property $A_2$:
\begin{theorem}\label{atwo}
 Each compact set $K$ of class $\mathfrak{C}_0$ has  property $A_2$.
 \end{theorem}
 This can be thought of as a
generalization of  Proposition~\ref{thesis}. The second result gives
examples of sets with property $A_3$:
\begin{theorem}\label{athree}
Each compact set $K$ of class $\mathfrak{C}_2$ has property $A_3$.
\end{theorem}
This again can be thought of as a generalization of
Proposition~\ref{forst}.

\subsection{}In \aref{goodpairsection} below,
we  prove a general result (Theorem~\ref{goodpairapprox} ) regarding
approximation of sections of submersions over a set $K$ which admits
a decomposition $K=K_1\cup K_2$ into a ``good pair" $(K_1,K_2)$ of
compact sets. This is a direct generalization of the results in
\cite{michiganpaper}, $\S$4.1 for maps into manifolds. Proofs of all
new statements are given in detail, although they are similar to the
proofs in the case of maps. This is the tool  for patching maps
which is used to prove Theorems~\ref{onedim} and \ref{athree}.

In  \aref{athreesection},  we give proofs of Lemma~\ref{zerodim},
Theorem~\ref{onedim}, and Theorem~\ref{athree}. For the last two,
the crucial ingredient is three successive applications of
Theorem~\ref{goodpairapprox} to produce a section which is
holomorphic in a neighborhood of the given compact $K$. For
Theorem~\ref{athree} we also require a result from
\cite{michiganpaper} (Theorem~\ref{arcs}). In \aref{atwosection} we
deduce Theorem~\ref{atwo} from Theorem~\ref{athree}.  This requires
a result from \cite{macgregor:interpolation} regarding conformal
mappings continuous to the boundary.

It should be noted that  the results proved here have analogs for
the approximation and interpolation of $\mathcal{A}^k$ maps, i.e.,
maps which are $\mathcal{C}^k$ and holomorphic in the interior. The
proofs are exactly the same. For notational simplicity we stick to
the case of maps which are only continuous to the boundary. Also,
this article is largely self contained, with the exception of the
proofs of Theorem~\ref{arcs} and Proposition~\ref{univalent}.

\section{Approximation on Good Pairs}\label{goodpairsection}
\subsection{Some definitions}
 We will call a set $K\subset \cx$ {\em nicely contractible} if
there is a homotopy $c:[0,1]\times K\rightarrow K$ with the
following properties:
\begin{enumerate}
\item for each $t$, the map $z\mapsto c(t,z)$ is in $\cts{K}{\cx}$,
\item $z\mapsto c(1,z)$ is the identity map on $K$.
\item there is a $z_0\in K$ such that $c(0,z)\equiv z_0$.
\item for $t\not=0$, $z\mapsto c(t,z)$ maps $\intr{K}$ into
$\intr{K}$.
\end{enumerate}
Of course, convex sets are nicely contractible, as are strongly
star-shaped sets (these are sets $K$ such that the maps $c(t,z)$ may
be taken as dilations with stretch factor $t$ and center $z_0\in
K$.) Another important class of nicely contractible sets are arcs.
The  property of nicely contractible sets which is used in the proof
of Lemma~\ref{cartan} (and therefore Theorem~\ref{goodpairapprox} is
the following: if $K$ is nicely contractible and $\m$ is a connected
complex manifold, then $\cts{K}{\m}$ is in fact contractible. We can
take the contraction to be the map $\phi_t$ from $\cts{K}{\m}$ to
itself given by $\phi_t(f)(z)= f(c(t,z)).$

Let $K_1$ and $K_2$ be compact subsets of $\cx$. We will say that
$(K_1,K_2)$ is a {\em good pair} if the following hold:
\begin{enumerate}
\item $\overline{K_1\setminus K_2}\cap \overline{K_2\setminus
K_1}=\emptyset$.
\item $K_{1,2}\assign K_1\cap K_2$ has finitely many connected components,
 each of which is nicely contractible.
\end{enumerate}

Let  $\m$ be a complex manifold of complex-dimension $n$, and
$\phi:\m\rightarrow\cx$ be a submersion. We will say that an open
set $U\subset \m$ is  {\em $\phi$-adapted}, if $U$ is biholomorphic
to an open set in $\cx^n$, and there are holomorphic coordinates
$(z_1,\cdots,z_n)$ on $U$ such that with respect to these
coordinates (and the standard coordinate on $\cx$), the map $\phi$
takes the form $(z_1,\cdots,z_n)\mapsto z_n$. (This is of course the
same as saying $z_n=\phi|_U$.)
\subsection{The main result}
We introduce the following notation and definitions:

\begin{enumerate}
\item Let $K=K_1\cup K_2$, where $(K_1,K_2)$ is a good pair.
\item Let $\m$ be a complex manifold and $\phi:\m\rightarrow\cx$ be a
holomorphic submersion such that $\phi(\m)\supset K$.

\item Let $B\subset\cx$ be compact and such that $B\cap
K_1=\emptyset$, and each  function $g\in\mathcal{A}(K_2)$ can be
approximated uniformly on $K_2$ by functions in $\mathcal{A}(K_2\cup
B)$ ;i.e., if $g\in \mathcal{A}(K_2)$ and $\epsilon>0$, then there
is a $g_\eta\in \mathcal{A}(K_2\cup B)$ such that
$\abs{g-g_\eta}<\eta$ on $K_2$.
\item Let $\mathcal{P}$ be a finite subset of $K$.
\end{enumerate}

We now state the following result regarding the approximation of
sections of $\phi$ over $K$:
\begin{theorem}\label{goodpairapprox}
With $K_1, K_2,\m,\phi$, $B$, and  $\mathcal{P}$  as above, let
$s\in\ctss{K}{\m}$ be a section of $\phi$ such that each $s(K_j)$
has a $\phi$-adapted neighborhood in $\m$ for $j=1,2$. Then, given
$\eta>0$, there is an $s_\eta\in \ctss{K}{\m}$ such that
$\dist(s,s_\eta)<\eta$ on $K$, $s_\eta $ extends as a holomorphic
section of $\phi$ to a neighborhood $B_\eta$ of $K_2\cap B$, and for
each $p\in\mathcal{P}$, we have $s_\eta(p)= s(p)$. Moreover,
$s_\eta(K_2\cup B_\eta)$ has a $\phi$-adapted neighborhood in $\m$.

\end{theorem}
The proof of \ref{goodpairapprox} will be reduced to the solution of
a non-linear patching problem in Euclidean space by the use of
coordinate in neighborhoods of $s(K_j)$. We now describe the
ingredients required in the proof . The most important is the
following version (resembling that A. Douady in
\cite{douady:modules}, pp.47-48)of H. Cartan's lemma on holomorphic
matrices.
\begin{lemma}\label{cartan} Let $\mathfrak{G}$ be a complex
connected Lie Group, and let $(K_1, K_2)$ be a good pair,  and let
$g\in \mathcal{A}(K_{1,2},\mathfrak{G})$. Then, for $j=1,2$ there
are $g_j\in \mathcal{A}(K_j,\mathfrak{G})$ such that $g=g_2\cdot
g_1$ on $K_{1,2}$.
\end{lemma}
For a proof, see \cite{michiganpaper}, Lemma~4.4, where it is
assumed that $\mathfrak{G}=GL_n(\cx)$, and each component of the
intersection is star shaped, but the proof is valid in general.

 If $\mathcal{P}\subset \cx$ is a finite set, we will let
$\mathcal{A}^\mathcal{P}(K,{\cx^n})$ denote the closed subspace of
the Banach space $\cts{K}{\cx^n}$ consisting of those functions
which vanish at each $p\in K\cap\mathcal{P}$. We will require the
following version of the solution of the additive Cousin problem
continuous to the boundary (the case with $\mathcal{P}=\emptyset$ is
in fact used in the proof of Lemma~\ref{cartan} ):
\begin{lemma}
\label{additive} Let $K_1,K_2$ be a compact subsets of the plane
such that $\overline{K_1\setminus K_2}\cap\overline{K_2\setminus
K_1}=\emptyset$, and let $\mathcal{P}$ be a finite subset of the
plane. There exist bounded linear maps $T_j:
\mathcal{A}^\mathcal{P}(K_{1,2},\cx)\rightarrow
\mathcal{A}^\mathcal{P}(K_j,\cx)$ such that for any function $f$ in
$\mathcal{A}^\mathcal{P}(K_{1,2},\cx)$ we have on $K_{1,2}$,
\begin{equation}
\label{additive_eq}
 T_1f +  T_2f = f,
\end{equation}
\end{lemma}
\begin{Proof} We reduce the problem to a $\dbar$
equation in the standard way. Let $\chi$ be a smooth cutoff which is
1 near $\overline{K_1\setminus K_2}$ and 0 near
$\overline{K_2\setminus K_1}$. Let $\lambda\assign f.\Db{\chi}$, so
that $\lambda\in \mathcal{A}(K_{1,2},\cx)$. Let
\[ \Lambda_f(z) =\frac{1}{2\pi i}\int_{K_{1,2}}\frac{\lambda(\zeta)}{\zeta-z}
d\overline{\zeta}\wedge d\zeta\]
Let $q_f$ be the Lagrange interpolation polynomial with the property
that $q_f(p)=\Lambda_f(p)$ for $p\in \mathcal{P}$. Observe that both
$\Lambda_f$ and $q_f$ are linear in $f$ and continuous in the sup
norm when restricted to compact sets.

 We can now define (assuming $(1-\chi).f=0$ where $\chi=1$ even if $f$ is not
 defined):
 \[ (T_1 f)(z) = (1-\chi(z)).f(z) +
 \Lambda_f(z) -q_f(z)\]

and (assuming $\chi.f=0$ where $\chi=0$ even if $f$ is not defined):
 \[ (T_2 f)(z) = \chi(z).f(z)
 - \Lambda_f(z) +q_f(z),\]
Since $T_jf$ is clearly holomorphic (resp. continuous) where $f$ is,
the result follows.
\end{Proof}
We will use the following standard result regarding Banach spaces,
which can be proved by iteration (see \cite{lang:realanalysis} pp.~
397-98):

\begin{lemma}
\label{surjective} In a metric space $X$, let $B_X(p,r)$ denote the
open ball in $X$ of radius $r$ centered at $p$.
  Let $\mathcal{E}$ and $\mathcal{F}$ be Banach Spaces and let
$\Phi: B_{\mathcal{E}}(p,r)\rightarrow\mathcal{F}$ be a
$\mathcal{C}^1$ map. Suppose there is a constant $C>0$ such that:
\begin{itemize}
\item
for each $h\in B_{\mathcal{E}}(p,r)$, the linear operator
 $\Phi'(h):\mathcal{E}\rightarrow\mathcal{F}$ is surjective and the
equation $\Phi'(h)u=g$ can be solved for $u$ in $\mathcal{E}$ for
all $g$ in $\mathcal{F}$ with the estimate
$\norm{u}_{\mathcal{E}}\leq C \norm{g}_{\mathcal{F}}$.
\item
for any $h_1$ and $h_2$ in $B_{\mathcal{E}}(p,r)$ we have $\norm{
\Phi'(h_1)-\Phi'(h_2)}\leq \frac{1}{2C}.$
\end{itemize}
Then,
\[ \Phi(B_{\mathcal{E}}(p,r)) \supset B_\mathcal{F}\left( \Phi(p),\frac{r}{2C}\right).\]
\end{lemma}

We will use the three lemmas above to give a proof of  the following
result (the main component of which is due to Rosay,
\cite{rosay:preprint}, and comments in \cite{rosay:korean})
regarding the solution of a non-linear Cousin problem (see Lemma 4.5
of \cite{michiganpaper}.) It is simply a translation of
Theorem~\ref{goodpairapprox} to coordinates.

\begin{lemma}\label{rosay}
Let $\omega$ be an open subset of $\cx^n$ and let
$\mathfrak{F}:\omega\rightarrow\cx^n$ be a holomorphic immersion.
Assume that $\mathfrak{F}$ preserves the last coordinate, i.e.,
$\mathfrak{F}$ is of the form $\mathfrak{F}(z_1,\cdots,z_n)=
\left(F(z_1,\cdots,z_n),z_n\right)$ for some map
$F:\omega\rightarrow\cx^{n-1}$.

Let $(K_1,K_2)$ denote a good pair of compact subsets of $\cx$, let
$\mathcal{P}$ be a finite subset of $K=K_1\cup K_2$, and suppose for
each $p\in K_2\setminus K_1$ we are given a point
$q(p)\in\cx^{n-1}$. Let $u_1\in \cts{K_1}{\cx^n}$ be such that
\begin{itemize}\item
$u_1(K_{1,2})\subset\omega$, and
\item  $u_1$ is of the form $u_1(z)=
\left(t_1(z),z\right)$ where $t_1:K_1\rightarrow\cx^{n-1}$
\end{itemize}
 Given any
$\epsilon>0$, there exists $\delta>0$ such that  if
$u_2\in\cts{K_2}{\cx^n}$ be such that
\begin{itemize}
\item $\norm{u_2-\mathfrak{F}\circ u_1}<\delta$ on $K_2\cap K_2$,
\item for $p\in \mathcal{P}\cap \left(K_{1,2}\right)$, we have
$u_2(p)=\mathfrak{F}(u_1(p))$,
\item for $p\in K_2\setminus K_1$ we have $u_2(p)=(q(p),p)\in\cx^n$, and
\item $u_2(z) =(t_2(z),z)$, where $t_2:K_2\rightarrow\cx^{n-1}$
\end{itemize}

then for $j=1,2$  there exist $v_j\in
\mathcal{A}^\mathcal{P}(K_j,\cx^n)$ such that
\begin{itemize}
\item $\norm{v_j}<\epsilon$,
\item $u_2+v_2=\mathfrak{F}(u_1+v_1)$, and
\item the last coordinate function of each of $v_1$ and $v_2$ is 0.
\end{itemize}
\end{lemma}

We now give a proof of Lemma~\ref{rosay}.

 \begin{Proof}
In order to apply Lemma~\ref{surjective} we choose the Banach spaces
$\mathcal{E}$, $\mathcal{F}$ and the map $\Phi$ as follows:
\begin{itemize}
\item  For $j=1,2$, let $\mathcal{B}_j$ be the closed
subspace of the Banach space $\mathcal{A}^\mathcal{P}(K_j,\cx^n)$
consisting of those maps whose last coordinate function is 0, i.e.,
$\mathcal{B}_j=\mathcal{A}^\mathcal{P}(K_j,\cx^n)\oplus 0$. We now
let $\mathcal{E}$ the Banach space
$\mathcal{B}_1\oplus\mathcal{B}_2$, which we endow with the norm
$\norm{\cdot}_{\mathcal{E}} \assign
\max\left(\norm{\cdot}_{\mathcal{A}(K_1,\cx^n)},\norm{\cdot}_{\mathcal{A}(K_2,\cx^n)}\right)$.
\item Let $\mathcal{F}$ be the Banach space $\mathcal{A}^\mathcal{P}(K_{1,2},\cx^{n-1})$.

\item Let $\pi:\cx^{n}\rightarrow\cx^{n-1}$ denote the projection on
the first $n-1$ coordinates, and let the open subset $\mathcal{U}$
of $\mathcal{E}$ be given by $\{(w_1,w_2): (u_1+w_1)(K_{1,2})\subset
\omega\}$. (Observe that $\{0\}\times\mathcal{B}(K_2,\cx^n)\subset
\mathcal{U}$.) Let $P=(P_1,\ldots,P_n)$  be  an $n$-tuple of
polynomials such that $P_n(z)\equiv z$, for $p\in K_{1,2}\cap
\mathcal{P}$, we have $P(p)=\mathfrak{F}(u_1(p))$, and for $p\in
K_2\setminus K_1$ we have $(P_1(p),\ldots, P_{n-1}(p))=q(p)$.

Let the map $\Phi:\mathcal{U}\rightarrow\mathcal{F}$ be given by
\[ \Phi(w_1,w_2)\assign \pi\circ\left[(P+w_2)|_{K_{1,2}}-
\mathfrak{F}\circ ((u_1+w_1)|_{K_{1,2}})\right].\] (Observe that,
since $\mathfrak{F}$ preserves the last coordinate, the last
coordinate function of $(P+w_2)|_{K_{1,2}}- \mathfrak{F}\circ
((u_1+w_1)|_{K_{1,2}})$ is actually 0. So the precomposition with
$\pi$ simply drops a coordinate which is identically $0$.)
\end{itemize}

 A computation shows that $\Phi'(w_1,w_2)$ is the bounded linear map from
$\mathcal{E}$ to $\mathcal{F}$ given by
\[(v_1,v_2)\mapsto \pi\circ \left[v_2|_{K_{1,2}}-
\mathfrak{F}'((u_1+w_1)|_{K_{1,2}})(v_1|_{K_{1,2}})\right].\]
Observe that $w_2$ plays no role whatsoever in this expression, and
therefore $\Phi'(w_1,w_2)\in BL(\mathcal{E},\mathcal{F})$ is in fact
a smooth function of $w_1$ alone, and we will henceforth denote it
by $\Phi'(w_1,*)$.

Let $\mathfrak{G}\subset GL_n(\cx)$ be the complex Lie subgroup  of
matrices of the form
\[ \left(\begin{array}{cc}A & b\\{\bf 0}&1\end{array}\right),\]
where $A\in GL_{n-1}(\cx)$, $b$ is an $1\times (n-1)$ column vector,
and ${\bf 0}$ is the zero row vector of size $(n-1)\times 1$. It is
easy to see that $\mathfrak{G}$ is connected.

 We  construct a right inverse to $\Phi'(u_1,*)$ .
Let $\gamma=\mathfrak{F}'\circ\left(u_1|_{K_{1,2}}\right)$. Since
$\mathfrak{F}$ preserves the last coordinate, $\gamma\in
\mathcal{A}(K_{1,2},\mathfrak{G})$, and thanks to Lemma~\ref{cartan}
above, we may write $\gamma= \gamma_2\cdot \gamma_1$, where
$\gamma_j\in \mathcal{A}(K_j,\mathfrak{G})$. (We henceforth suppress
the restriction signs.) Let
$j:\mathcal{F}\rightarrow\cts{K_{1,2}}{\cx^n}$ be the continuous
inclusion induced by the map $\cx^{n-1}\rightarrow\cx^{n}$ given by
$(z_1,\cdots, z_{n-1})\mapsto (z_1,\cdots,z_{n-1},0)$.
 For $g\in \mathcal{F}$,
 let \[ S(g)=(-\gamma_1^{-1} T_1(\gamma_2^{-1}j(g)),\gamma_2
T_2(\gamma_2^{-1}j(g))),\] where the $T_j$ are as in
equation~\ref{additive_eq} above (and extended to $\cx^n$
componentwise). Observe that $-\gamma_1^{-1} T_1(\gamma_2^{-1}j(g))$
and $\gamma_2 T_2(\gamma_2^{-1}j(g))$ vanish at each
$p\in\mathcal{P}$ and they belong to $\mathcal{B}_1$ and
$\mathcal{B}_2$ respectively. It is easy to see that
$S:\mathcal{F}\rightarrow\mathcal{E}$ is a bounded linear operator,
and a computation shows that $\Phi'(u_1,*)\circ S=\id_\mathcal{F}$.
Choose $\theta>0$ so small so that if $w_1\in
B_\mathcal{F}(u_1,\theta),$ \begin{enumerate}
\item$\norm{\Phi'(w_1,*)-
\Phi'(u_1,*)}_{\tmop{op}}< \frac{1}{8\norm{S}}$ (possible by
continuity), and
 \item  the
equation $\Phi'(w_1,*)u=g$ can be solved with the estimate
$\norm{u}\leq 2 \norm{S} \norm{g}$, (possible from  the fact that
small perturbations of  a surjective linear operator is still
surjective.) \end{enumerate}
 Consequently, if $\epsilon<\theta$ and
$\tilde{u}_2\in\mathcal{B}(K_2,\cx^n)$, for the ball
$B_{\mathcal{E}}((0,\tilde{u}_2),\epsilon)$ the hypothesis of
Lemma~\ref{surjective} are verified with $C=2\norm{S}$. We have
therefore,
\begin{eqnarray*}
\Phi\left(B_\mathcal{E}((0,\tilde{u}_2),\epsilon)\right)&\supset&
 B_\mathcal{F}\left(\Phi(0,\tilde{u}_2),\frac{\epsilon}{2C}\right)\\
& =&  B_\mathcal{F}\left(P+\tilde{u}_2-
\mathfrak{F}({u_1}),\frac{\epsilon}{2C}\right).
\end{eqnarray*}
So, if $\norm{P+\tilde{u}_2- F(u_1)}< \frac{\epsilon}{4C}$, we have
$0\in\Phi\left(B_\mathcal{E}((0,\tilde{u}_2),\epsilon)\right)$. This
is exactly the conclusion required, since any $u_2$ such that
$u_2(p)=P(p)$ for each $p\in K_2\cap \mathcal{P}$  can be written as
$u_2= P+\tilde{u}_2$ for some
$\tilde{u}_2\in\mathcal{B}(K_2,\cx^n)$.
 \end{Proof}
 We will now prove Theorem~\ref{goodpairapprox}.
\begin{Proof}
We omit the restriction signs on maps for notational clarity. For
$j=1,2$ let $V_j$ be a $\phi$-adapted neighborhood of $s(K_j)$, and
let $\mathfrak{F}_j:V_j\rightarrow\cx^n$ be a coordinate system such
that $\mathfrak{F}_j(z)=\left(F_j(z),\phi(z)\right)$. Let
$\mathfrak{F}=\mathfrak{F}_2\circ\mathfrak{F}_1^{-1}$ be the
associated transition function.Then $\mathfrak{F}$ is a
biholomorphism from the open set $\omega=\mathfrak{F}_1(V_1\cap
V_2)$ onto the open set $\mathfrak{F}_2(V_2\cap V_1)$, and
 $\mathfrak{F}$ preserves the last coordinate, i.e.,
$\mathfrak{F}$ is of the form $\mathfrak{F}(z_1,\cdots,z_n)=
\left(F(z_1,\cdots,z_n),z_n\right)$ for some map
$F:\omega\rightarrow\cx^{n-1}$. Any section of $\phi$ over $K_j$ is
represented in the coordinate system $\mathfrak{F}_j$  by a map of
the form $t_j:K_j\rightarrow\cx^n$, where $t_j(z)=
(\tilde{t}_j(z),z)$, with $\tilde{t}_j$ a map from $K_j$ to
$\cx^{n-1}$. Also, for $j=1,2$, given maps $t_j:K_j\rightarrow\cx^n$
of the form $t_j=(\tilde{t}_j,z)$,  they    glue together to form a
section over $K_1\cup K_2$ (i.e. there is a section $\lambda$ of
$\phi$ over $K_1\cup K_2$ such that $t_j =\mathfrak{F}_j\circ
\lambda$) iff $t_2 = \mathfrak{F}\circ t_1$.

 Since $B\cap
K_1=\emptyset$ by hypothesis, the pair of compact sets $(K_1,
K_2\cup B)$ is good. Let $\epsilon>0$, and let
$u_1=\mathfrak{F}_1\circ s$. Observe that $u_1(K_{1,2})\subset
\omega$, and $u_1$ is of the form $u_1(z)= (t_1(z),z)$, where
$t_1:K_1\rightarrow\cx^{n-1}$. Then Lemma~\ref{rosay} gives a
$\delta>0$ corresponding to $u_1$, the good pair $(K_1, K_2\cup B)$,
$\omega$ and $\mathfrak{F}$.

Let $w_2 = \mathfrak{F}_2\circ (s|_{K_2})$, then
$w_2\in\cts{K_2}{\cx^n}$, and is of the form $w_2(z)=
(\tilde{w}_2(z),z)$. Thanks to the hypothesis regarding uniform
approximation of functions in $\cts{K_2}{\cx}$ by functions in
$\cts{K_2\cup B}{\cx}$, we can find a $u_2\in\cts{K_2\cup B}{\cx^n}$
of the same form as $w_2$ such that $\norm{u_2 -
\mathfrak{F}(u_1)}<\delta$ on $K_{1,2}$ (Since $\mathfrak{F}(u_1)=
w_2$ on $K_{1,2}$.) Further we may assume that for
$p\in\mathcal{P}\cap K_2$, we have $u_2(p)=w_2(p)$ and that the last
coordinate function of $u_2$ is $z$.

Then, by Proposition~\ref{rosay},  there is a
$v_1\in\mathcal{A}^\mathcal{P}(K_1,\cx^n)$ and a
$v_2\in\mathcal{A}^\mathcal{P}(K_2\cup B,\cx^n)$, such that
$\norm{v_j}<\epsilon$ and $u_2 + v_2 = \mathfrak{F}(u_1 + v_1)$, and
the last coordinate functions of the $v_j$ are  0's. Hence the maps
$u_1+v_1$ and $u_2+v_2$ glue together to form a section of $\phi$
(which we call $\tilde{s}_\epsilon$ given by
\[
\tilde{s}_{\epsilon}\assign \left\{ \begin{array}{ccc}
\mathfrak{F}_1^{-1}(u_1+v_1)& \tmop{ on } & K_1\\
\mathfrak{F}_2^{-1}(u_2+v_2) &  \tmop{ on } & K_2 \tmop{  and ~near
} K_2\cap B
\end{array} \right.
\]
Clearly, $\tilde{s}_\epsilon$ is in $\ctss{K_1\cup K_2}{\m}$, and
extends to a holomorphic map near $K_2\cap B$. Moreover,
$\tmop{dist}(\tilde{s}_\epsilon, s)= O(\epsilon)$. By construction,
we have $s_\epsilon(p)=s(p)$. Therefore, given $\eta>0$, we can find
$s_\eta$ with required properties.
\end{Proof}
\section{Sets with property $A_3$}\label{athreesection}
In this section we prove Lemma~\ref{zerodim}, and
Theorems~\ref{onedim} and \ref{athree}. In each we show that a
certain set $K$ has property $A_3$. We will let $\m$ be a complex
manifold, $\phi:\m\rightarrow\cx$ a holomorphic submersion such that
$\phi(\m)\supset K$, $\mathcal{P}$ a finite subset of $K$, and $s$ a
section of $\phi$, $s\in\ctss{K}{\m}$. We let $\epsilon>0$, and want
to show that there is an $s_\epsilon\in\hols{K}{\m}$ such that
$\dist(s,s_\epsilon)<\epsilon$ and $s_\epsilon(p)=s(p)$ for each
$p\in\mathcal{P}$.
\subsection{Proof of Lemma~\ref{zerodim}}
If $K$ is not finite, it can be written as a disjoint union of
finitely many singletons and a closed subset $C$ homeomorphic to the
Cantor middle-third set. (See e.g., \cite{ems:gtopology1}, pp.
108-109.) In particular, $\cx\setminus K$ is connected, and $K$ has
property $A_1$.

Observe that $s\in\ctss{K}{\m}$  is simply a continuous section of
$\phi$ over $K$. We can cover $K$ by finitely many sets $\{U_j\}$
such that each $s(\overline{U_j})$ has a $\phi$-adapted neighborhood
in $\m$. We can choose a refinement $\{V_j\}$ of this cover, such
that the $V_j$'s are pairwise disjoint. Observe that $K_j\assign
K\cap \overline{V_j}$ has property $A_1$. Now, with respect to the
$\phi$-adapted coordinate system around $s(K_j)$, the map $s$ has a
representation on $K_j$ of the form $s(z)= (\tilde{s}(z),z)$, where
$\tilde{s}$ is continuous and takes values in $\cx^{n-1}$.
Approximating $\tilde{s}$ by a holomorphic map in a neighborhood of
$K_j$ our result follows.

\subsection{Proof of Theorem~\ref{onedim}}
We first introduce some combinatorial preliminaries. Let $vw$ be an
edge in a graph $\Gamma$. We can construct a new graph $\Gamma'$
with one more vertex by ``splitting the edge $vw$." More formally,
if $V$ is the vertex set of $\Gamma$, and $E$ its edge set, the
vertex set of $\Gamma'$ is $V\cup\{u\}$ (where $u\not\in V$), and
the edge set is $\left(E\setminus\{ vw\}\right)\cup \{vu,uw\}$.

Recall that a graph is {\em $n$-colorable}, if there is an
assignment of $n$ colors to its vertices so that no two adjacent
vertices have the same color. By a classical result of K\"{o}nig
(see \cite{harary:book}, Theorem 12.1), the condition that a graph
is $2$-colorable is that it does not have a circuit (closed non-self
intersecting path) of odd length. We now have the following:
\begin{lemma}\label{bicolorable}
Given any finite graph $G$, there are edges $e_1,\ldots, e_k$ such
that the graph $G'$ obtained after successively splitting these
edges is  2-colorable.
\end{lemma}
\begin{Proof} Let $V$ be the vector space over the field $\mathbb{Z}/2\mathbb{Z}$
with basis the set of edges of $G$. Given a circuit $C$ in $G$,
traversing the edges $x_1,x_2,\ldots, x_m$, associate with it an
element $c$ of $V$ given by $c= x_1+x_2+\cdots+x_m$. Let $Z$ be the
subspace of $V$ spanned by all  $c$ for each circuit $C$. Pick a
basis $c_1,\ldots, c_k$. Since any circuit can be written as a
linear combination of the $c_j$'s it is sufficient to split any edge
of those $c_j$ which have an odd number of summands. Moreover, from
the fact that the $c_j$ are linearly independent over
$\mathbb{Z}/2\mathbb{Z}$ it follows that each $c_j$ has an edge not
contained in any  $c_k$, $k\not=j$. The result follows.
\end{Proof}
Now we can prove Theorem~\ref{onedim}. Note that as in the proof of
Lemma~\ref{onedim} above, $\intr{K}=\emptyset$, so  $s$ is simply a
continuous section of $\phi$ over $K$.

For $z\in K$, there is a neighborhood $U_z$ of $z$  in $K$ such that
$s(\overline{U_z})$ has a $\phi$-adapted neighborhood in $\m$.
Select a finite subcover $\mathcal{F}$ of $\{U_z\}_{z\in K}$, and
let $\mathcal{G}$ be a refinement of $\mathcal{F}$ which has only
double intersections (this is possible since $\tmop{dim}(K)=1$, and
$K$ is connected.) We can write $\mathcal{G}=\{V_1,\ldots,V_M\}$.

We associate to  $\mathcal{G}$ a graph $N_\mathcal{G}$ by taking the
nerve: the vertices $V_i$ and $V_j$ are joined by an edge in
$N_\mathcal{G}$ iff $V_i\cap V_j\not=\emptyset$. Let $V$ and $W$ be
open sets in $\mathcal{G}$ such that $V\cap W\not=\emptyset$. We
define a new open cover $S_{VW}(\mathcal{G})$ of $K$ in the
following way. Replace the sets $V$ and $W$ by $V'$, $U$ and $W'$,
where $V'\subset V$, $W'\subset W$, and $U$ is a neighborhood of
$\overline{V\cap W}$ such that $V'\cup U\cup W'= V\cup W$, and
$V\cap U\not=\emptyset$, $U\cap W\not=\emptyset$ but $V\cap W=
\emptyset$. If $U$ is chosen small enough, then
$\left(\mathcal{G}\setminus\{V,W\}\right)\cup \{V',U,W'\}$ is again
an open cover of $K$ with only double intersections. This is our
$\mathcal{G}'=S_{VW}(\mathcal{G})$. At the level of nerves
$N_{\mathcal{G}'}$ is obtained by splitting  the edge $VW$ in the
graph $N_{\mathcal{G}}$, as in Lemma~\ref{bicolorable}. Thanks to
Lemma~\ref{bicolorable} we can obtain after finitely many rounds of
edge-splitting  a  new cover $\mathcal{H}$, with only double
intersections such that the corresponding nerve $N_\mathcal{H}$ is
$2$-colorable.  It also follows that for each $U\in \mathcal{H}$,
the set $s(\overline{U})$ has a $\phi$-adapted neighborhood in $\m$.
Further, by shrinking the $U$'s we can assume that whenever an
intersection $U\cap U'$ of two sets in $\mathcal{H}$ is non-empty,
it is contractible. Let us call the colors used in coloring
$N_\mathcal{H}$ red and blue. We can now define
\[ K_1=\bigcup_{{\stackrel{U\in\mathcal{H}}{ U\mbox{~{\scriptsize{red}}}}}}
\overline{U},\mbox{~~and~~~}
K_2=\bigcup_{\stackrel{U\in\mathcal{H}}{U\mbox{~{\scriptsize{blue}}}}}
\overline{U}.\]
It is easy to verify that $(K_1,K_2)$ is a good pair, and for
$j=1,2$ the set $s(K_j)$ has a $\phi$-adapted neighborhood in $\m$.

We also define sets $K_1'\subset K_1$ and $K_2'\subset K_2$ in the
following way. For $U,U'\in\mathcal{H}$,  let $\mathcal{W}_{UU'}$ be
a simply connected neighborhood of $U\cap U'$ in $\cx$ with
$\mathcal{C}^2$ boundary. Further we can assume that
$\overline{\mathcal{W}_{UU'}\cap K}$ is contractible and
$s(\overline{\mathcal{W}_{UU'}\cap K})$ is contained in a
$\phi$-adapted open set of $\m$. We set
\[\mathcal{W}= \bigcup_{\stackrel{U,U'\in\mathcal{H}}{U\cap
U'\not=\emptyset}}\mathcal{W}_{UU'},\] so that $\mathcal{W}$ is a
neighborhood in $\cx$ of the points in $K$ which are contained in
two sets of the cover $\mathcal{H}$.

Let $V\in\mathcal{H}$. We set
\[V^s= V\setminus \overline{W} \]

So that $V^s\subset V$, and points of $V^s$ do not belong to any
other set of $\mathcal{H}$ apart from $V$. We set:
\[ K_1'=\bigcup_{{\stackrel{U\in\mathcal{H}}{ U\mbox{~{\scriptsize{red}}}}}}
\overline{U^s},\mbox{~~and~~~}
K_2'=\bigcup_{\stackrel{U\in\mathcal{H}}{U\mbox{~{\scriptsize{blue}}}}}
\overline{U^s}.\]

We now apply Theorem~\ref{goodpairapprox} to the good pair $(K_1,
K_2)$, and obtain an
 $s_1\in \hols{K\cup\mathcal{B}_1}{\m}$, where
$\mathcal{B}_1$ is a neighborhood of $K_1'$ in $\cx$, such that
$\dist(s,s_1)<\frac{\epsilon}{3}$, and $s(p)=s_1(p)$ for
$p\in\mathcal{P}$. We can further assume that $\mathcal{B}_1$ is a
disjoint union  of simply connected neighborhoods of the sets $V^s$
for $V$ red,  $\partial\mathcal{B}_1$ is smooth, and
$\partial\mathcal{B}_1$ and $\partial\mathcal{W}$ meet transversely
at each point of intersection.

Observe now that $(K_1\cup \mathcal{B}_1, K_2)$ is a good pair, and
we can  apply Theorem~\ref{goodpairapprox} to it. We obtain a
$s_2\in \ctss{K\cup \mathcal{B}_1\cup \mathcal{B}_1}{\m}$ with
$\dist(s_2,s_1)<\frac{\epsilon}{3}$, and $s_2(p)= s_1(p)$ for
$p\in\mathcal{P}$. As before, $\mathcal{B}_2$ is a disjoint union of
simply connected neighborhoods of the sets $V^s$ for $V$ blue,
$\partial\mathcal{B}_2$ is smooth, and $\partial\mathcal{B}_2$ and
$\partial\mathcal{W}$ meet transversely at each point of
intersection.

We set $\mathcal{B}=\mathcal{B}_1\cup \mathcal{B}_2$ so that
$\partial\mathcal{B}=\partial\mathcal{B}_1\cup
\partial\mathcal{B}_2$ meets $\partial\mathcal{W}$ transversely, so
that $\partial\mathcal{B}\cap\partial\mathcal{W}$ is a finite set.
Let $\mathcal{U}$ be a simply connected neighborhood of $K$ in $\cx$
such that all points of $\partial\mathcal{B}\cap\partial\mathcal{W}$
lie outside $\mathcal{U}$ and $\partial\mathcal{U}$ meets both
$\partial\mathcal{B}$ and $\partial\mathcal{W}$ transversely. Note
that $K_1\cup K_2\subset\mathcal{W}\setminus\mathcal{B}$. Let
\[ L_1 = \partial\mathcal{U}\cap\partial\mathcal{B},\]
and
\[ L_2 = \partial\mathcal{U}\cap \partial\mathcal{W}.\]
It is easy to verify that $(L_1,L_2)$ is a good pair, and each of
$s(L_j)$ has a $\phi$-adapted neighborhood in $\m$. Setting
$L=L_1\cup L_2$, we apply Theorem~\ref{goodpairapprox} to
$(L_1,L_2)$ to obtain a map $s_3\in\ctss{L\cup \mathcal{C}}{\m}$,
where  $\mathcal{C}$ is a neighborhood of $K_1\cap K_2$ in $\cx$,
$\dist(s_3,s_2)<\frac{\epsilon}{3}$, and $s_3(p)=s_2(p)$ for all
$p\in\mathcal{P}$. Clearly, $s_3\in \hols{K}{\m}$, and we are done.

\subsection{Proof of Theorem~\ref{athree}}
In \aref{arcsection} we record a result that will be used in the
proof. In \aref{smoothathree} we use Theorem~\ref{arcs} and
Theorem~\ref{goodpairapprox} to give a proof of Theorem~\ref{athree}
in the special case when $N=1$, i.e. $K$ is the closure of a
smoothly bounded domain. As in the proof of Theorem~\ref{onedim}
this uses a ``Triple bumping." Finally, in \aref{athreeproofsection}
we complete the proof by reducing the general case to the case
considered in \aref{smoothathree}.
\subsubsection{Arcs in Complex Manifolds}
\label{arcsection} If $X$ is a differentiable manifold, and $\alpha$
is an arc (continuous injective map from $[0,1]$) which is at least
$\mathcal{C}^1$, we will say $\alpha$ is {\em embedded} if for each
$t\in[0,1]$, we have $\alpha'(t)\not=0$. A proof of the following
result, required in the proof of Theorem~\ref{athree}, can be found
in \cite{michiganpaper}, Theorem 2:

\begin{theorem}
\label{arcs}

Let $\m$ be a complex manifold, and $\alpha$ be an arc in $\m$.
Assume that
\begin{enumerate}
\item there is a complex-valued submersion $\phi$ defined in a
neighborhood of $\alpha$ in $\m$ such that $\phi\circ\alpha$ is a
$\mathcal{C}^1$ embedded arc in $\cx$.
\item  there is a finite subset $P\subset [0,1]$ such that $\alpha$ is
$\mathcal{C}^3$ on $[0,1]\setminus P$.
\end{enumerate}
Then $\alpha$ has a $\phi$-adapted neighborhood in $\m$.
\end{theorem}

\subsubsection{Case of a  smooth domain}\label{smoothathree} In this
section we prove the following special case of Theorem~\ref{athree}
(Compare with Proposition~\ref{forst}):
\begin{proposition}\label{smoothathreeprop}
Let $\Omega\Subset\cx$ be a domain with $\mathcal{C}^2$ boundary.
Then $\overline{\Omega}$ has property $A_3$. \end{proposition} The
first step in the proof is the following application of Sard's
Theorem:
\begin{lemma}\label{sard1}
 There is a unit vector $\mathbf{v}$ in the plane such that:
\begin{itemize}\item every straight line in the plane parallel to
$\mathbf{v}$ meets $\partial\Omega$ in only finitely many points.
\item the number of straight lines parallel to $\mathbf{v}$ which
are tangent to $\partial\Omega$ is finite.
\end{itemize}
In fact these hold for almost all unit vectors $\mathbf{v}$ in the
unit circle $S^1$ (with respect to the standard measure on $S^1$.)
\end{lemma}
\begin{Proof}Fix an arbitrary orientation on $\partial\Omega$, and define
the Gauss map $G:\partial\Omega\rightarrow S^1$ by mapping the point
$z\in\partial\Omega$ to the unit tangent vector $G(z)$ to
$\partial\Omega$ at the point $z$. This is a $\mathcal{C}^1$ map,
and the set of its critical values is of measure 0 by Sard's
Theorem. If $\mathbf{v}$ is any regular value of $G$, such that
$-\mathbf{v}$ is also a regular value,  it follows that the sets
$G^{-1}(\mathbf{v})$ and $G^{-1}(-\mathbf{v})$ are discrete in
$\partial\Omega$ (since $G$ is a diffeomorphism near each of them.)
Since $\partial\Omega$ is compact, it follows that the sets
$G^{-1}(\pm\mathbf{v})$   are finite. From this, the conclusions
follow immediately.
\end{Proof}

After a rotation if required, we will assume that $\mathbf{v}$ is
vertical, i.e., $\mathbf{v}=\pm i$. Now we introduce some notation.
for $c\in\rl$, denote by $L(c)$ the vertical straight line
$\Re(z)=c$ in $\cx$, and for $a<b$, denote by $M[a,b]$ the vertical
strip $\{z\in\cx\colon a\leq\Re(z)\leq b\}$. Also set
$l(c)=L(c)\cap\overline{\Omega}$ and $m[a,b]=
M[a,b]\cap\overline{\Omega}$. We will assume without loss of
generality that $\Omega\subset M[0,1]$.

Thanks to Lemma~\ref{sard1} and the choice $\mathbf{v}=\pm i$, it
follows that for each $c$, the set $l(c)$ has only finitely many
components, each of which is either a point, or a line segment.
Also, $l(c)\cap\partial\Omega$ is finite.  We now make the following
observation.

\begin{lemma}\label{delta}
There is a $\delta>0$ such that for each $c\in[0,1]$, the set
$s(m[c-\delta,c+\delta])$ has a $\phi$-adapted neighborhood in $\m$.
\end{lemma}
\begin{Proof} First we show that for each $c\in[0,1]$, the set
$s(l(c))$ has a $\phi$-adapted neighborhood in $\m$. Since $s$ is
injective (it has a left-inverse $\phi$), it is sufficient to show
that each component of $l(c)$ has a $\phi$-adapted neighborhood in
$\m$. For a component of $l(c)$ that reduces to a point, this is
trivial. Therefore, consider a component which is a (vertical) line
segment. We can think of this component as an arc in the plane, and
parameterize it as $\lambda(t) = z_0+ iat$, where $z_0$ is the lower
end point of the segment, and $a$ is a positive real number. Then
$s\circ\lambda$ is a continuous arc in $\m$, which is real analytic
except at finitely many points, and $\phi$ is a submersion from $\m$
to $\cx$ such that $\phi\circ\left(s\circ\lambda\right)=\lambda$.
Thanks to Proposition~\ref{arcs} above, $s(\lambda)$ has a
$\phi$-adapted neighborhood in $\m$. Therefore, $s(l(c))$ has a
$\phi$-adapted neighborhood in $\m$.

It follows that there is a $\delta_c$ such that
$s(m[c-\delta_c,c+\delta_c])$ has a $\phi$-adapted neighborhood in
$\m$. The uniform choice of $\delta$ follows by compactness.
\end{Proof}

We will also require the following simple fact (a proof may be found
in \cite{michiganpaper} , Observation 4.8 )
\begin{lemma}\label{strip}
Let $u$ and $v$ be real valued  $\mathcal{C}^1$ functions defined on
a neighborhood of $0$ in $\rl$ such that for each $x$, we have
$u(x)<0<v(x)$. Then there is an $\eta>0$ such that for
$0<\theta\leq\eta$, the vertical strip
\[ S \assign  \{ (x,y)\in \rl^2\colon x\in [-\theta,\theta],
 {u}(x)\leq y \leq {v}(x) \} \]
is star shaped with respect to the origin.
\end{lemma}

The proof of Proposition~\ref{smoothathree} will parallel that of
Proposition~\ref{onedim} in that both require a ``Triple bumping",
i.e. three successive applications of Theorem~\ref{goodpairapprox}.
However, the good pairs are obtained by different methods.

Let $\mathcal{E}\subset[0,1]$ be the set of $c$ such that the line
$L(c)$ is tangent to some component of $\partial\Omega$.
$\mathcal{E}$ is finite by Lemma~\ref{sard1}. If
$c\not\in\mathcal{E}$,  $L(c)$ meets $\partial\Omega$ transversely
at each point of intersection, so that (1) each component of $l(c)$
is a line segment, and (2) (by Lemma~\ref{strip} ) there is a
$\theta_c$ such that for $\epsilon\leq\theta_c$, each component of
$m([c-\epsilon,c+\epsilon])$ is strongly star shaped.

By compactness, we can find $c_1<c_2<\ldots<c_M$, and $\eta_j>0$,
$j=1,\ldots,M$, such that if $m_j= m[c_j-\eta_j,c_j+\eta_j]$, we
have $\overline{\Omega}=\bigcup_{j=1}^M m_j$. We can assume that the
$\eta_j<\frac{\delta}{100}$, where $\delta$ is as in
Lemma~\ref{delta}. We see that for each $j$, $s(m_j)$ has a
$\phi$-adapted neighborhood in $\m$. We will further impose the
following conditions on the $m_j$'s
\begin{enumerate}
\item Each point of $\overline{\Omega}$ is in at most two of the
$m_j$'s (this can be done by shrinking the $\eta_j$'s.) Therefore
$m_j\cap m_k=\emptyset$ if $\abs{j-k}>1$.
\item Each $c\in\mathcal{E}$ occurs in the list $\{c_j\}_{j=1}^M$.
\item Each component of $m_j\cap m_{j+1}$ is strongly star shaped.
Observe that by the previous step, $m_j\cap m_{j+1}$ does not
contain any $l(c)$ for $c\in\mathcal{E}$. Therefore, this can be
achieved by shrinking the $\eta_j$'s.
\end{enumerate}
Now let
\[ K_1=\bigcup_{j\mbox{~~{\scriptsize{odd}}}}m_j ,\mbox{~~and~~}
 K_2=\bigcup_{j\mbox{~~{\scriptsize{even}}}}m_j.\]
 It is easy to see that $(K_1,K_2)$ is a good pair. Since each
 $s(m_j)$ has a $\phi$-adapted neighborhood, and $K_1,K_2$ are
 disjoint union of $m_j$'s, it follows that each of $s(K_1)$ and
 $s(K_2)$ has a $\phi$-adapted neighborhood in $\m$.

 Let $I_1\subset \partial\Omega\cap K_1$ and $I_2\subset
 \partial\Omega\cap K_2$ be such that (1) $ I_1\cap K_2= I_2\cap K_2
 =\emptyset,$ and (2) Each connected component
 $\partial\Omega\setminus (I_1\cup I_2)$ is
 contained in a vertical strip of width $\frac{\delta}{2}$, where
 $\delta$ is as in Lemma~\ref{delta}. (This is possible, since
 $\eta_j<\frac{\delta}{100}$.)

Let $B$ be a neighborhood of $I_1$ such that $B\cap K_2=\emptyset$.
We apply Theorem~\ref{goodpairapprox} to the good pair $(K_1,K_2)$
to obtain an $s_1\in \ctss{\overline{\Omega}\cup B_1}{\m}$, where
$B_1$ is a neighborhood of $I_1$ contained in $B$, such that
$\dist(s,s_1)<\frac{\epsilon}{3}$, and $s(p)=s_1(p)$ for $p\in
\mathcal{P}$. Observe that for $B_1$ small enough, $s_1(K_1\cup
B_1)$ is contained in a $\phi$-adapted open set of $\m$, and
$(K_1\cup B_1, K_2)$ is again a good pair. We now apply
Theorem~\ref{goodpairapprox} again to this good pair to obtain an
$s_2\in\ctss{\overline{\Omega}\cup B_1\cup B_2}{\m}$ where $B_2$ is
a neighborhood of $I_2$, such that
$\dist(s_1,s_2)<\frac{\epsilon}{3}$, and $s_1(p)= s_2(p)$ for
$p\in\mathcal{P}$.

Let $I_3=\partial\Omega\setminus (I_1\cup I_2).$ By construction,
each connected component of $I_3$ is contained in a vertical strip
of width $\frac{\delta}{2}$. Let $\Omega'$ be a domain with
$\mathcal{C}^2$ boundary such that $\overline{\Omega}\cup B_1\cup
B_2\supset \Omega'\supset \Omega\cup I_1\cup I_2$ (i.e. $\Omega'$ is
obtained by smoothly bumping $\Omega$ along $I_1$ and $I_2$.) We can
assume that $\Omega'$ and $\Omega$ are so close that for $\delta$ as
in Lemma~\ref{delta},  and any $c$, the set
$s(M[c-\delta,c+\delta]\cap\overline{\Omega'})$ is contained in a
$\phi$-adapted open set in $\m$. We can now repeat the constructions
that gave us $K_1$ and $K_2$ to obtain a good pair $(L_1,L_2)$ such
that (1) $L_1\cup L_2=\overline{\Omega'}$, (2) each of $s(L_1)$ and
$s(L_2)$ has a $\phi$-adapted neighborhood in $\m$, and (3)
$I_3\subset L_1\setminus L_2$. We can now apply
Theorem~\ref{goodpairapprox} to the good pair $(L_1,L_2)$ to obtain
an $s_3\in\ctss{\overline{\Omega'}\cup C}{\m}$, where $C$ is a
neighborhood of $I_3$ in $\cx$, such that
$\dist(s_3,s_2)<\frac{\epsilon}{3}$ and $s_3(p)=s_2(p)$, for each
$p\in\mathcal{P}$. Then $s_3\in\hols{\overline{\Omega}}{\m}$, and we
are done.

\subsubsection{End of proof of Theorem~\ref{athree}}
\label{athreeproofsection} Let $K$ be of class
$\mathfrak{C}_2$. Recall that $K=\cup_{i=1}^N\overline{\Omega_i}$,
and for $i\not =j$ , we have $\partial\Omega_i$ and
$\partial\Omega_j$ meet at a set of finitely many points $P_{ij}$.
Set $P=\cup_{i\not=j}P_{ij}$. We can refer to the points in $P$ as
{\em nonsmooth points} of $\partial K$.

 The proof in this section is very
similar to that in \aref{smoothathree}. The first step is to
establish the following version of Lemma~\ref{sard1} for this case:
\begin{lemma}\label{sard2}
Let $K$ be of class $\mathfrak{C}_2$. Then there is a unit vector
$\bf{v}$ in the plane with the following properties:
\begin{enumerate}\item
every straight line in the plane parallel to $\mathbf{v}$ meets
$\partial K$ in only finitely many points.
\item The number of straight lines parallel to $v$ which are tangent to
$\partial K$ at smooth points is finite.
\item Let $p\in P_{ij}$ be a non-smooth point of $\partial K$. Then,
the straight line through $p$ parallel to $\mathbf{v}$ is transverse
to both $\partial\Omega_i$ and $\partial\Omega_j$ at $p$.
\end{enumerate}
\end{lemma} \begin{Proof}Applying lemma~\ref{sard1} separately to
each $\partial\Omega_j$, we conclude that properties (1) and (2)
hold for almost all unit vectors $\mathbf{v}$. Let $p$ be a
non-smooth point of $\partial K$, so that for some $i,j$, we have
$p\in P_{ij}=\partial\Omega_i\cap
\partial\Omega_j$. Let $\mathbf(t)(p)$ be a common unit tangent
vector to $\partial\Omega_i$ and  $\partial\Omega_j$ at the point
$p$. We can choose $\mathbf{v}\not= \pm \mathbf{t}(p)$ for all $p\in
P_{ij}$ for all $i$ and $j$.
\end{Proof}

As in \aref{smoothathree} we can assume that $\mathbf{v}=\pm i$. Let
$L(c)$ and $M[a,b]$ have the same meaning as in the last section,
and set $l'(c)= L(c)\cap K$, $m'[a,b]= M$. We can assume that
$K\subset M[0,1]$.

We let $\mathcal{E}$ be the finite set of points $c\in[0,1]$ such
that either (1) $L(c)$ is tangent to $\partial\Omega_i$ for some
$i$, or (2) $L(c)$ passes through a nonsmooth point of $\partial K$.
Arguing as in the previous section, we can find $c_1<c_2<\ldots c_M$
and $\eta_j>0$, $j=1,\ldots, M$, such that if
$m_j'=m'[c_j-\eta_j,c_j+\eta_j]$, we have $K=\cap_{j=1}^M m_j'$ and
each $s(m_j'$ has a $\phi$-adapted neighborhood in $\m$. We can
further impose the following conditions, (the first three are just
as in the last section, the last is a  new condition):
\begin{enumerate}
\item Each point in $K$ is contained in at most two of the $m_j'$'s,
i.e.,  $m_j'\cap m_k'=\emptyset$ if $\abs{j-k}>1$.
\item Each $c\in\mathcal{E}$ occurs among the $\{c_j\}_{j=1}^M$.
\item Each component of $m_j'\cap m_{j+1}'$ is strongly star shaped
(Note that, thanks to the last step, in this case we have ensured
that $m_j'\cap m_{j+1}'$ does not contain any nonsmooth points.)
\item Each nonsmooth point is contained in an $m_j'$ with an even
$j$. This can be ensured by introducing additional $c_j's$. Observe
that if a nonsmooth point $p\in m_j'$, then $p\not\in m_{j-1}'$ and
$p\not \in m_{j+1}'$.
\end{enumerate}

As in the proof of Proposition~\ref{smoothathreeprop}
\[ K_1=\bigcup_{j\mbox{~~{\scriptsize{odd}}}}m_j' ,\mbox{~~and~~}
 K_2=\bigcup_{j\mbox{~~{\scriptsize{even}}}}m_j'.\]
 It is easy to see that $(K_1,K_2)$ is a good pair. Since each
 $s(m_j')$ has a $\phi$-adapted neighborhood, and $K_1,K_2$ are
 disjoint union of $m_j'$'s, it follows that each of $s(K_1)$ and
 $s(K_2)$ has a $\phi$-adapted neighborhood in $\m$. Moreover, the
 set $P$ of nonsmooth points of $\partial K$ is contained in
 $K_2\setminus K_1$.

 Let $B$ be a neighborhood of the nonsmooth points
 $P$ such that $B\cap K_1=\emptyset$.
 Thanks to Theorem~\ref{goodpairapprox} we can find an
 $s'\in\ctss{K\cup B'}{\m}$ (where $B'$ is a neighborhood of $P$
 contained in $B$), such that $\dist(s,s')<\frac{\epsilon}{2}$, and
 $s'(p)=s(p)$ for each $p\in\mathcal{P}$.

Now we can find an open set $\Omega'$ with $\mathcal{C}^2$ boundary
such that $\overline{\Omega'}\subset K\cup B'$ but $\Omega'\supset
\Omega\cup P$. Then $s'\in\cts{\overline{\Omega'}}{\m}$, and thanks
to Proposition~\ref{smoothathreeprop}, we can find $s_\epsilon\in
\hols{\overline{\Omega'}}{\m}$ such that
$\dist(s_\epsilon,s')<\frac{\epsilon}{2}$ and $s_\epsilon(p)=s'(p)$
for $p\in \mathcal{P}$. Since $s_\epsilon\in \hols{K}{\m}$, the
proof is complete.

\section{Proof of Theorem~\ref{atwo}}\label{atwosection}
\subsection{$\mathcal{A}$-equivalence}
Let $K$ and $K'$ be compact subsets of $\cx$. We will say that $K$
and $K'$ are {\em $\mathcal{A}$-equivalent} if there is a
homeomorphism $\chi:K\rightarrow K'$ such that $\chi|_{\intr{K}}$ is
a conformal map of $\intr{K}$ onto $\intr{K'}$. We will call  $\chi$
an {\em $\mathcal{A}$-equivalence } from $K$ to $K'$. A well-known
example of $\mathcal{A}$-equivalence is the following
(\cite{tsuji:bible}, Theorems IX.35 and IX.2):
\begin{lemma}\label{circular} Let $\Omega$ be a {\em Jordan Domain},
Then there is a domain $\omega$ in the plane bounded by {\em
circles} such that $\overline{\Omega}$ and $\overline{\omega}$ are
$\mathcal{A}$-equivalent.
\end{lemma}
The significance of this notion in the current investigation is
explained by the following observation:
\begin{lemma}\label{aeq} Suppose that two compact sets $K$ and $L$ in $\cx$ are
$\mathcal{A}$-equivalent. If $K$ has property $A_2$ and $L$ has
property $A_1$ then $L$ has property $A_2$.
\end{lemma}
\begin{Proof}
Let $\m$ be a complex manifold, $f\in\cts{L}{\m}$, and $\mathcal{P}$
be a finite subset of $K$. We want to approximate $f$ by maps $f_n$
in $\hol{L}{\m}$ such that $f_n(p)=f(p)$ for $p\in\mathcal{P}$. Let
$\chi:K\rightarrow L$ be an $\mathcal{A}$-equivalence,let
$\mathcal{Q}=\chi^{-1}(\mathcal{P})$, and let $g=f\circ\chi$. Then
$g\in\cts{K}{\m}$, and consequently there is a sequence
$g_n\in\hol{K}{\m}$ such that $g_n\rightarrow g$ uniformly, and
$g_n(q)=g(q)$ for $q\in\mathcal{Q}$. Let $\zeta=\chi^{-1}$, so that
$\zeta\in\cts{L}{\cx}$. Since $L$ has property $A_1$,  we can find
$\zeta_n\in\hol{L}{\cx}$ such that $\zeta_n\rightarrow\zeta$ on $L$,
with $\zeta_n(p)=\zeta(p)$ for $p\in\mathcal{P}$. Then $f_n\assign
g_n\circ\zeta_n\in\hol{L}{\m}$,  $f_n\rightarrow f$ uniformly, and
$f_n(p)=f(p)$ for $p\in\mathcal{P}$.
\end{Proof}
\subsection{Proof of Theorem~\ref{atwo}} Thanks to Lemma~\ref{aeq}
and Theorem~\ref{athree}, it is sufficient to prove the following
result:
 \begin{theorem}\label{touching}
For each $K$ in $\mathfrak{C}_0$ there is an $L$ in
$\mathfrak{C}_\omega$ such that $K$ and $L$ are
$\mathcal{A}$-equivalent.
\end{theorem}

We will require two results, the first from Combinatorics.  A vertex
$v$ of a graph $G$ is said to be a {\em cutpoint} of $G$ if the
graph $G^{\{v\}}$ obtained by removing from $G$ the vertex $v$ along
with all edges incident at $v$ has at least one more connected
component than $G$ has. (So for example, if $G$ is connected,
$G^{\{v\}}$ is {\em not} connected.)  We will need the following
elementary  fact.
\begin{lemma}\label{graph}(\cite{harary:book}, Theorem~3.4, p. 29)
Let $G$ be a graph with more than one vertex. Then there are at
least two vertices of $G$ which are {\em not} cutpoints.
\end{lemma}
The second result is the following boundary interpolation theorem
for conformal maps, due to MacGregor and Tepper
(\cite{macgregor:interpolation}, Theorem~1). $\Delta\subset\cx$ is
the open unit disc, and $\{z_1,z_2,\ldots,z_n\}\subset
\partial\Delta$ and $\{w_1,w_2,\ldots,w_n\}\subset \partial\Delta$
are given finite subsets of the unit circle.
\begin{proposition}\label{univalent}  There is a function $f$ which is analytic and univalent in the
union of $\Delta$ and a neighborhood of $\{z_1,z_2,\ldots,z_n\}$ and
continuous on $\overline{\Delta}$ such that $f(z_k)=w_k$ for
$k=1,\ldots,n$. Furthermore, $\abs{f(z)}=1$ if $\abs{z}=1$ and $z$
is sufficiently near any of the points $z_k$, and also
$f(\Delta)\subset\Delta.$
\end{proposition}
For a compact connected set $K$ in the plane, by the outer boundary
we mean the boundary of the unbounded component of $\cx\setminus K$
(this is also a component of $\partial K$.)
 We use Proposition~\ref{univalent} to prove the following lemma.
\begin{lemma}\label{interpolation}
Let $\Omega\Subset\cx$ be a Jordan domain, and let $\gamma$ be its
outer boundary.  Suppose we are given a finite set of points
$\{z_1,\ldots, z_n\}$ on $\gamma$ and the same number of points
$\{w_1,\ldots,w_n\}$ on the unit circle $\partial\Delta$. Then there
is a continuous map
$f:\overline{\Omega}\rightarrow\overline{\Delta}$ such that
\begin{enumerate}
\item $f|_\Omega$ is conformal,
\item $f(z_k)= w_k$, for $k=1,\ldots, n$,
\item let $W=f(\Omega)$. Then $\partial W$ is $\mathcal{C}^\omega$,
and at each $w_k$, $\partial W$ is tangent to $\partial \Delta$.
\end{enumerate}
\end{lemma}
\begin{Proof}
By Lemma~\ref{circular} there is a domain $D$ bounded by circles and
an $\mathcal{A}$-equivalence
$f_0:\overline{\Omega}\rightarrow\overline{D}$. After applying an
inversion of the plane if required, we can assume further than the
outer boundary $\gamma$ of $\Omega$ is mapped onto the outer
boundary of $D$, which we may assume is the unit circle
$\partial\Delta$. Set $z_k'= f_0(z_k)$ for $k=1,\ldots,n$.

Let $D'$ be a simply connected open set in $\cx$ with
$\mathcal{C}^\infty$ boundary such that $\Delta\subset D'$,
$\partial D'\cap \partial \Delta =\{z_1',\ldots, z_n'\}$ where at
each $z_k'$, the boundaries $\partial D'$ and $\partial \Delta$ are
tangent to each other. Let $f_1:D'\rightarrow \Delta$ be a conformal
map of $D'$ onto $\Delta$. Since $\partial D'$ has
$\mathcal{C}^\infty$ boundary, $f_0$  extends to a diffeomorphism of
the closures. Set $z_k''=f_1(z_k')$, and $D''=f_1(D')$. Observe that
$f_1\circ f_0$ maps $\Omega$  to a subdomain $\Omega''$ of $\Delta$
and the boundary $\partial\Omega''$ is tangent to $\partial\Delta$
at each point of intersection $z_k''$.

We now apply Proposition~\ref{univalent} to obtain a continuous
$f_2:\overline{\Delta}\rightarrow\overline{\Delta}$ such that
$f_2(z_k'')=w_k$, $f_2$ is conformal on the union of $\Delta$ with a
neighborhood of $\{z_1'',\ldots,z_n''\}$, and $f_2$ maps a piece of
$\partial\Delta$ near each $z_k''$ onto a piece of $\partial\Delta$
near $w_k$. It follows immediately that if $W=f_2(\Omega'')$,
$\partial W$ meets $\partial\Delta$ tangentially at each $w_k$. We
set $f\assign f_2\circ f_1\circ f_0$. The properties claimed are
easily verified.

\end{Proof}
 We now prove Theorem~\ref{touching}.

\begin{Proof} It is clear that we only need to consider the case in
which $K$ is connected. We use induction on $N$, the number of
summands of $K$.

When $N=1$, the result is reduced to  Lemma~\ref{circular}. Now
suppose that the result has been proved for some $N\geq 1$, and let
$K= \cup_{i=1}^{N+1} \overline{\Omega_j}$. Let $G$ be a graph whose
vertices $v_i$ correspond to the sets $\overline{\Omega_i}$, and
there is an edge connecting $v_i$ and $v_j$ iff
$\overline{\Omega_i}\cap \overline{\Omega_j}\not=\emptyset$. Since
we have assumed that $K$ is connected, it follows that $G$ is a
connected graph. Thanks to Lemma~\ref{graph} above, we can assume
(after a renumbering of the vertices of $G$) that $v_{N+1}$ is {\em
not} a cutpoint of $G$. Let $K'=\cup_{i=1}^N \overline{\Omega_j}$,
and let $P\subset K'$ be the finite set
$K'\cap\overline{\Omega_{N+1}}=K'\cap\partial\Omega_{N+1}$.  Then
$K'$ is connected, therefore  is contained in exactly one connected
component $U$ of $\cx\setminus\Omega_{N+1}$. Let $\gamma=\partial
U$. Clearly $\gamma$ is a connected component of
$\partial\Omega_{N+1}$. It follows that the set $P\subset\gamma$.
Moreover, as $\overline{\Omega_{N+1}}$ is connected, it follows that
$\Omega_{N+1}$ is contained in exactly one component of
$\cx\setminus K'$.

We claim that {\em we can assume that $\gamma$ is the outer boundary
of $\Omega_{N+1}$}.  To show this it is sufficient to show that for
some $\mathcal{A}$-equivalence $\Phi:K\rightarrow
\hat{K}\subset\cx$, $\Phi(\gamma)$ is the outer boundary of
$\Phi(\Omega_{N+1})$.

 If $\gamma$ not already the outer boundary of $\Omega_{N+1}$
 let $z_0\in U\setminus K'$, where $U\Subset\cx$ is the component of
$\cx\setminus \Omega_{N+1}$ which contains $K'$ (then
$\gamma=\partial U$). Let $\rho>0$ be small enough so that
$B_\cx(z_0,\rho)\Subset U\setminus K'$, and define the inversion
$\Phi:\cx\setminus\{z_0\}\rightarrow\cx$ by
\[ \Phi(z) = \frac{\rho^2}{z-z_0}.\]
Then $\Phi(K)$ is contained in the ball $B_\cx(0,\rho)$, and since
$z_0$ is mapped to the point at infinity, it follows that
$U\setminus K'$ is mapped to the unbounded component of
$\cx\setminus K$. Since $\gamma=\partial U$, we see that
$\gamma\subset \partial \left( U\setminus K'\right)$, so that
$\Phi(\gamma)$ is the outer boundary of $\Phi(\Omega_{N+1})$.

Now, by induction hypothesis, there is an
$L'\in\mathfrak{C}_\omega$, and an $\mathcal{A}$-equivalence
$\chi':K'\rightarrow L'$. Using Lemma~\ref{univalent} we will extend
$\chi'$ to an $\mathcal{A}$-equivalence $\chi$ defined on  $K=K'\cup
\overline{\Omega_{N+1}}$.

Let us write $P=\{\zeta_1,\ldots, \zeta_n\}$, and let $\zeta_k' =
\chi'(\zeta_k)$. Then the $\zeta_k'$'s lie at the boundary of a
single connected component $U$ of $\cx\setminus L'$. Let $U'$ be a
simply connected domain, $U\subset U'$ such that $\partial U'$ is
$\mathcal{C}^\omega$ and passes through each $\zeta_k'\in\partial
U$, and further at each $\zeta_k'$, $\partial U'$ is tangent to
$\partial U$, i.e. to $\partial L$. Let $\theta$ be a conformal map
of $U'$ onto the disc $\Delta$. Then $\theta$ extends to a
holomorphic map of a neighborhood of $\overline{U'}$. We set
$w_k=\theta(\zeta_k')$.

Thanks to Lemma~\ref{interpolation} above, there is a map
$\lambda\in\cts{\overline{\Omega_{N+1}}}{\cx}$ such that
$\lambda|_{\Omega_{N+1}}$ is conformal,
$\lambda(\Omega_{N+1})\subset \Delta$ $f(\zeta_k)=w_k$,
$\partial(\lambda(\Omega_{N+1}))$ is real analytic and tangent to
$\partial\Delta$ at each point $w_k$. We can now define
\[ \chi\assign\left\{\begin{array}{ccc} \chi'&\mbox{on}& K'\\
\theta^{-1}\circ\lambda &\mbox{on} &
\overline{\Omega_{N+1}}.\end{array}\right.
\]
Let $\omega_{N+1}=\chi(\Omega_{N+1})$, and $L = L'\cup
\overline{\omega_{N+1}}$. Then $L$ is in $\mathfrak{C}_\omega$, and
$\chi$ is an $\mathcal{A}$-equivalence between $K$ and $L$.
\end{Proof}
\section{A Problem} We conclude this article by
stating an open problem. An  solution  will lead to a clearer
picture of sets with properties  $A_2$ and $A_3$.

{\em Is it possible to prove an analog of
Theorem~\ref{goodpairapprox} for {\em three } sets $K_1$, $K_2$,
$K_3$?} That is, given $s: K\rightarrow\m$, (where $K=K_1\cup
K_2\cup K_3$), such that $s(K_j)$ lies in a subset of $\m$
homeomorphic to an open set in $\cx^n$, obtain an approximation to
$s$, after assuming reasonable hypotheses. In particular, we should
have $K_1\cap K_2\cap K_3\not=\emptyset$.

Such a result will be necessary if  we want to avoid the use of
results like Theorem~\ref{arcs} two prove approximation results for
two-dimensional sets. Observe that the use of Theorem~\ref{arcs}
resulted in assumptions regarding the smoothness of the sets on
which we want to do approximation.

\section{Acknowledgements} It is a great pleasure to thank
Andr\'{e} Boivin and Rasul Shafikov for all their help in the
preparation of this article.


\begin{thebibliography}{10}

\bibitem{ems:gtopology1}
A.V. Arkhangel'ski\v{i} and V.V. Fedorchuk.
\newblock {\em General Topology I}, volume~17 of {\em Encyclopaedia of
  Mathematical Sciences}.
\newblock Springer-Verlag, 1988.

\bibitem{michiganpaper}
Debraj Chakrabarti.
\newblock Coordinate neighborhoods of arcs and the approximation of maps into
  (almost) complex manifolds.
\newblock {\em Michigan Math. J.}, 55, 2007.
\newblock To Appear, available at {\tt{arxiv.org}}, math.CV/0605496.

\bibitem{demailly:nonsteinbundle}
Jean-Pierre Demailly.
\newblock Un exemple de fibr\'e holomorphe non de {S}tein \`a fibre {${\bf
  C}\sp{2}$} ayant pour base le disque ou le plan.
\newblock {\em Invent. Math.}, 48(3):293--302, 1978.

\bibitem{douady:modules}
Adrien Douady.
\newblock Le probl\`eme des modules pour les sous-espaces analytiques compacts
  d'un espace analytique donn\'e.
\newblock {\em Ann. Inst. Fourier (Grenoble)}, 16(fasc. 1):1--95, 1966.

\bibitem{bddforstneric:approximation}
B.~Drinovec-Drnovsek and F.~Forstneric.
\newblock Approximation of holomorphic mappings on strongly pseudoconvex
  domains.
\newblock Available online at {\tt{arxiv.org}},math.CV/0607185.

\bibitem{bddforstneric:holocurves}
B.~Drinovec-Drnovsek and F.~Forstneric.
\newblock Holomorphic curves in complex spaces.
\newblock {\em Duke Math. J.}, To Appear.
\newblock available at {\tt{arxiv.org}}, math.CV/0604118.

\bibitem{forstneric:noncritical}
Franc Forstneri{\v{c}}.
\newblock Noncritical holomorphic functions on {S}tein manifolds.
\newblock {\em Acta Math.}, 191(2):143--189, 2003.

\bibitem{harary:book}
Frank Harary.
\newblock {\em Graph theory}.
\newblock Addison-Wesley Publishing Co., Reading, Mass.-Menlo Park,
  Calif.-London, 1969.

\bibitem{henkinleiterer:book}
G.~M. Henkin and J.~Leiterer.
\newblock {\em Theory of functions on complex manifolds}, volume~60 of {\em
  Mathematische Lehrb\"ucher und Monographien, II. Abteilung: Mathematische
  Monographien [Mathematical Textbooks and Monographs, Part II: Mathematical
  Monographs]}.
\newblock Akademie-Verlag, Berlin, 1984.

\bibitem{lang:realanalysis}
Serge Lang.
\newblock {\em Real and functional analysis}, volume 142 of {\em Graduate Texts
  in Mathematics}.
\newblock Springer-Verlag, New York, third edition, 1993.

\bibitem{macgregor:interpolation}
T.~H. MacGregor and D.~E. Tepper.
\newblock Finite boundary interpolation by univalent functions.
\newblock {\em J. Approx. Theory}, 52(3):315--321, 1988.

\bibitem{rosay:preprint}
Jean-Pierre Rosay.
\newblock Approximation of non-holomorphic maps, and {P}oletsky theory of
  discs.
\newblock Unpublished Preprint.

\bibitem{rosay:korean}
Jean-Pierre Rosay.
\newblock Approximation of non-holomorphic maps, and {P}oletsky theory of
  discs.
\newblock {\em J. Korean Math. Soc.}, 40(3):423--434, 2003.

\bibitem{tsuji:bible}
M.~Tsuji.
\newblock {\em Potential theory in modern function theory}.
\newblock Maruzen Co. Ltd., Tokyo, 1959.

\bibitem{vitushkin}
A.~G. Vitu{\v{s}}kin.
\newblock Conditions on a set which are necessary and sufficient in order that
  any continuous function, analytic at its interior points, admit uniform
  approximation by rational fractions.
\newblock {\em Dokl. Akad. Nauk SSSR}, 171:1255--1258, 1966.
\newblock English translation in Soviet Mat. Dokl. {\em 7}(1966), 1622-1625.

\end{thebibliography}
\end{document}